\pdfoutput=1
\documentclass[letterpaper, 10 pt, conference]{ieeeconf}  

\IEEEoverridecommandlockouts                              

\usepackage{mathrsfs}
\usepackage[font={small}]{caption}
\usepackage{scalerel}
\usepackage{url}
\usepackage{bbold}
\usepackage{amsfonts}
\usepackage{amsmath,amssymb,amsfonts}
\usepackage{graphicx}
\usepackage{wrapfig}
\usepackage{mathrsfs} 
\usepackage{algorithm,algorithmic}
\usepackage{array}
\usepackage{times}
\usepackage{url}
\usepackage{cite}
\newcommand{\citet}[1]{\cite{#1}}
\usepackage{upgreek}
\usepackage{float}
\usepackage{longtable}
\usepackage{color}
\usepackage{wasysym}
\usepackage{grffile}
\usepackage{stackengine}

\allowdisplaybreaks

\newcommand{\0}{\mathbb{0}}

\newcommand{\ub}{\mathbf{u}}
\newcommand{\wb}{\mathbf{w}}
\newcommand{\xb}{\mathbf{x}}

\newcommand{\Ab}{\mathbf{A}}
\newcommand{\Bb}{\mathbf{B}}

\newcommand{\Hb}{\mathbf{H}}
\newcommand{\Ib}{\mathbf{I}}

\newcommand{\Kb}{\mathbf{K}}
\newcommand{\Lb}{\mathbf{L}}

\newcommand{\Pb}{\mathbf{P}}
\newcommand{\Qb}{\mathbf{Q}}
\newcommand{\Rb}{\mathbf{R}}

\newcommand{\Vb}{\mathbf{V}}

\newcommand{\Wb}{\mathbf{W}}
\newcommand{\Xb}{\mathbf{X}}

\newcommand{\Fc}{\mathcal{F}}

\newcommand{\Nc}{\mathcal{N}}

\newcommand{\Sc}{\mathcal{S}}

\newcommand{\argmin}{\text{argmin}}

\newcommand{\minimize}{\text{minimize}}

\newcommand{\norm}[1]{\left\lVert#1\right\rVert}

\newtheorem{theorem}{Theorem}

\newtheorem{definition}{Definition}

\newtheorem{lemma}[theorem]{Lemma}

\title{\LARGE \bf Distributed Online Linear Quadratic Control for \\ Linear Time-invariant Systems}

\author{Ting-Jui Chang and Shahin Shahrampour, {\it Senior Member}, {\it IEEE}  
\thanks{T.J. Chang and S. Shahrampour are with Wm Michael Barnes '64 Department of Industrial and Systems Engineering, Texas A\&M University, College Station, TX 77843, USA. 
        {\tt\small email:\{tingjui.chang,shahin\}@tamu.edu}.}%
}

\begin{document}

\maketitle
\thispagestyle{empty}
\pagestyle{empty}

\begin{abstract}
Classical linear quadratic (LQ) control centers around linear time-invariant (LTI) systems, where the control-state pairs introduce a quadratic cost with time-invariant parameters. Recent advancement in online optimization and control has provided novel tools to study LQ problems that are robust to time-varying cost parameters. Inspired by this line of research, we study the distributed online LQ problem for identical LTI systems. Consider a multi-agent network where each agent is modeled as an LTI system. The LTI systems are associated with decoupled, time-varying quadratic costs that are revealed sequentially. The goal of the network is to make the control sequence of all agents competitive to that of the best centralized policy in hindsight, captured by the notion of regret. We develop a distributed variant of the online LQ algorithm, which runs distributed online gradient descent with a projection to a semi-definite programming (SDP) to generate controllers. We establish a regret bound scaling as the square root of the finite time-horizon, implying that agents reach consensus as time grows. We further provide numerical experiments verifying our theoretical result. 
\end{abstract}

\section{Introduction}
In recent years, there has been a significant interest on problems arising at the interface of control and machine learning. Modern statistical and optimization algorithms have opened new avenues to rethink classical control problems, where linear quadratic (LQ) control (\cite{4309169, bertsekas1995dynamic, zhou1996robust}) is a prominent point in case. In its classical form, LQ control centers around LTI systems, where the control-state pairs introduce a quadratic cost with {\it time-invariant} parameters. For the infinite-horizon problem, the optimal controller has a closed-form solution, and it can be derived by solving the algebraic Riccati equation.

Fueled by applications in practical control problems, {\it online} LQ control has received a great deal of attention \citet{cohen2018online}. In this scenario, the environment is subject to unpredictable dynamics, making the cost functions {\it time-varying} in an arbitrary fashion. Examples include variable-supply electricity production and building climate control with time-varying energy costs. Motivated by the online optimization literature, online LQ casts the time-varying problem as a regret minimization, where the performance of the online algorithm is compared with that of the best fixed control policy in hindsight.

In this paper, we address the distributed online LQ problem for a network of identical LTI systems. Each system is modeled as an agent in a multi-agent network, associated with a decoupled, time-varying quadratic cost. The cost sequence for each agent can be chosen in an adversarial fashion and the agent observes the sequence on-the-fly. The goal of the network is to make the control sequence of all agents competitive to that of the best centralized policy in hindsight, captured by the notion of {\it regret}. We develop a distributed variant of the online LQ algorithm. At each iteration, agents run distributed online gradient descent \cite{yan2012distributed} to maintain an ideal steady-state covariance matrix. To do so, they need to perform a projection to an SDP and extract a feasible policy to generate the controllers. We prove that the individual regret can be bounded by $O(\sqrt{T})$, where $T$ is the total number of iterations. This implies that the agents reach consensus and collectively compete with the best fixed controller in hindsight. We finally provide simulation results verifying this theoretical property.

\subsection{Related Work}
\noindent
\textbf{Distributed LQ Control:}
Distributed linear quadratic regulator (LQR) has been widely studied in the control literature. Some works focus on multi-agent systems with known, identical decoupled dynamics. In \cite{4626964}, a distributed control design is proposed by solving a single LQR problem whose size matches the maximum vertex degree of the underlying graph plus one. The authors of \citet{6862471} derive the necessary condition for the optimal distributed controller, resulting in a non-convex optimization problem. The work of \cite{5299181} addresses a multi-agent network, where the dynamics of each agent is a single integrator. The authors of \cite{5299181} show that the computation of the optimal controller requires the knowledge of the graph and the initial information of all agents. Given the difficulty of precisely solving the optimal distributed controller, Jiao et al. \citet{8736845} provide the sufficient conditions to obtain sub-optimal controllers. All of the aforementioned works need global information such as network topology to compute the controllers. On the other hand, Jiao et al. \citet{8734804} propose a decentralized way to compute the controllers and show that the system will reach consensus. For the case of unknown dynamics, Alemzadeh et al. \cite{alemzadeh2019distributed} propose a distributed Q-learning algorithm for dynamically decoupled systems. There are other works focusing on distributed control without assuming identical decoupled sub-systems. Fattahi et al. \cite{fattahi2019efficient} study distributed controllers for unknown and sparse LTI systems. Furieri et al. \cite{furieri2020learning} address model-free methods for distributed LQ problems and provide sample-complexity bounds for problems with local gradient dominance property (e.g., quadratically-invariant problems). The work of \cite{furieri2020first} investigates the convergence of distributed controllers to a global minimum for quadratically invariant problems with first-order methods. 

\vspace{0.2cm}
\noindent
{\bf Classical LQ with Unknown Dynamics:} There is a recent line of research dealing with LQ control problems with unknown dynamics. Several techniques are proposed using (i) gradient estimation (see e.g., \cite{fazel2018global, malik2019derivative, 9130755, 9147749}) (ii) the estimation of dynamics matrices and derivation of the controller by considering the estimation uncertainty \cite{dean2019sample}, and (iii) wave-filtering \cite{hazan2017learning, arora2018towards}.

\vspace{0.2cm}
\noindent
\textbf{Online LQ Control:}
Recently, there has been a significant interest in studying linear dynamical systems with time-varying cost functions, where online learning techniques are applied. This literature investigates two scenarios: 
\begin{enumerate}
    \item \textbf{Known Systems:} As mentioned before, Cohen et al. \citet{cohen2018online} study the SDP relaxation for online LQ control and establish a regret bound of $O(\sqrt{T})$ for known LTI systems with time-varying quadratic costs. Agarwal et al. \citet{agarwal2019online} propose the disturbance-action policy parameterization and reduce the online control problem to online convex optimization with memory. They show that for adversarial disturbances and arbitrary  time-varying convex functions, the regret is $O(\sqrt{T})$. Agarwal et al. \cite{agarwal2019logarithmic} consider the case of time-varying strongly-convex functions and improve the regret bound to $O(\text{poly}(\text{log}T))$. Simchowitz et al. \citet{simchowitz2020improper} further extend the $O(\text{poly}(\text{log}T))$ regret bound to partially observable systems with semi-adversarial disturbances.
    \item \textbf{Unknown Systems:} For fully observable systems, Hazan et al. \citet{hazan2020nonstochastic} derive the regret of $O(T^{2/3})$ for time-varying convex functions with adversarial noises. For partially observable systems, the work of \cite{simchowitz2020improper} addresses the cases of (i) convex functions with adversarial noises and (ii) strongly-convex functions with semi-adversarial noises, and provide regret bounds of $O(T^{2/3})$ and $O(\sqrt{T})$, respectively. Lale et al. \citet{lale2020logarithmic} establish an $O(\text{poly}(\text{log}T))$ regret bound for the case of stochastic perturbations, time-varying strongly-convex functions, and partially observed states.
\end{enumerate}

Our work lies precisely at the interface of distributed LQR and online LQ, addressing distributed online LQ.  

\section{Problem Formulation}
\subsection{Notation} 
We use the following notation in this work:

{\small
\begin{tabular}{|c||l|}
    \hline
    $[n]$ & The set of $\{1,2,\ldots,n\}$ for any integer $n$ \\
    \hline
    $\text{Tr}(\cdot)$ & The trace operator\\
    \hline
    $\norm{\cdot}$ & Euclidean (spectral) norm of a vector (matrix)\\
    \hline
    $\mathrm{E}[\cdot]$ & The expectation operator\\
    \hline
    $[\Ab]_{ij}$ & The entry in the $i$-th row and $j$-th column of $\Ab$\\
    \hline
    \rule{0pt}{11pt}$\Ab\bullet\Bb$ & $\text{Tr}(\Ab^\top\Bb)$\\
    \hline
    $\Ab\succeq\Bb$ & $(\Ab-\Bb)$ is positive semi-definite\\
    \hline
\end{tabular}}

\subsection{Distributed Online LQ Control}
We consider a multi-agent network of $m$ identical LTI systems. The dynamics of agent $i$ is given as, 
$$\xb_{i,t+1} = \Ab\xb_{i,t}+\Bb\ub_{i,t}+\wb_{i,t},\quad i\in [m]$$
where $\xb_{i,t}\in\mathrm{R}^d$ and $\ub_{i,t}\in\mathrm{R}^k$ represent agent $i$'s state  and control (or action) at time $t$, respectively. Furthermore, $\Ab\in\mathrm{R}^{d\times d}$, $\Bb\in\mathrm{R}^{d\times k}$, and $\wb_{i,t}$ is a Gaussian noise with zero mean and covariance $\Wb\succ 0$. The noise sequence $\{\wb_{i,t}\}$ is independent over time and agents. 

Departing from the classical LQ control, we consider the {\it online} distributed LQ problem in this work. At round $t$, agent $i$ receives the state $\xb_{i,t}$ and applies the action $\ub_{i,t}$. Then, positive definite cost matrices $\Qb_{i,t}$ and $\Rb_{i,t}$ are revealed, and the agent incurs the cost $\xb_{i,t}^\top\Qb_{i,t}\xb_{i,t} + \ub_{i,t}^\top\Rb_{i,t}\ub_{i,t}$. Throughout this paper, we assume that $\text{Tr}(\Qb_{i,t}),\text{Tr}(\Rb_{i,t})\leq C$ for all $i,t$ and some $C>0$. Agent $i$ follows a policy that selects the control $\ub_{i,t}$ based on the observed cost matrices $\Qb_{i,1},\ldots,\Qb_{i,t-1}$ and $\Rb_{i,1},\ldots,\Rb_{i,t-1}$, as well as the information received from the agents in its neighborhood.

\vspace{.2cm}
\noindent
{\bf Centralized Benchmark:} In order to gauge the performance of any distributed LQ algorithm, we require a centralized benchmark. In this work, we focus on the {\it finite-horizon} problem, where for a centralized policy $\pi$, the cost after $T$ steps is given as 
\begin{align}\label{eq:cost-central}
    J_T(\pi)=\mathrm{E}\left[\sum_{t=1}^T {\xb_t^{\pi}}^\top\Qb_t\xb_t^{\pi} + {\ub_t^{\pi}}^\top\Rb_t\ub_t^{\pi}\right],
\end{align}
where $\Qb_t = \sum_{i=1}^m \Qb_{i,t}$ and $\Rb_t = \sum_{i=1}^m \Rb_{i,t}$, and the expectation is over the possible randomness of the policy as well as the noise. The superscript $\pi$ in $\ub_t^{\pi}$ and $\xb_t^{\pi}$ alludes that the state-control pairs are chosen by the policy $\pi$, given full access to cost matrices of all agents. Notice that in the {\it infinite-horizon} version of the problem with time-invariant cost matrices $(\Qb,\Rb)$, where the goal is to minimize $\lim_{T\rightarrow\infty}J_T(\pi)/T$, it is well-known that for a controllable LTI system $(\Ab,\Bb)$, the optimal policy is given by the constant linear feedback, i.e., $\ub_t^{\pi}=\Kb\xb_t^{\pi}$ for a matrix $\Kb\in \mathrm{R}^{k\times d}$.

\vspace{.2cm}
\noindent
{\bf Regret Definition:} The goal of a distributed LQ algorithm $\mathcal{A}$ is to mimic the performance of an ideal centralized algorithm using only {\it local} information. More formally, each agent $j$ locally generates the control sequence $\{\ub_{j,t}\}_{t=1}^T$, that is competitive to the best policy among a benchmark policy class $\Pi$. This can be formulated as minimizing the individual regret, which is defined as follows
\begin{align}\label{eq:regret}
    \text{Regret}_T^j(\mathcal{A})=J_T^j(\mathcal{A})-\min_{\pi\in\Pi}J_T(\pi),
\end{align}
for agent $j\in [m]$, where 
\begin{align*}
    J_T^j(\mathcal{A}) &= \mathrm{E}\left[\sum_{t=1}^T\sum_{i=1}^m {\xb_{j,t}^{\mathcal{A}}}^\top\Qb_{i,t}\xb_{j,t}^{\mathcal{A}} + {\ub_{j,t}^{\mathcal{A}}}^\top\Rb_{i,t}\ub_{j,t}^{\mathcal{A}}\right]\\
    &=\mathrm{E}\left[\sum_{t=1}^T {\xb_{j,t}^{\mathcal{A}}}^\top\Qb_{t}\xb_{j,t}^{\mathcal{A}} + {\ub_{j,t}^{\mathcal{A}}}^\top\Rb_{t}\ub_{j,t}^{\mathcal{A}}\right].
\end{align*}
A successful distributed algorithm is one that keeps the regret sublinear with respect to $T$. Of course, this also depends on the choice of the benchmark policy class $\Pi$, which is assumed to be the set of strongly stable policies (to be defined precisely in Section \ref{sec:stability}).

\vspace{.2cm}
\noindent
{\bf Network Structure:} The agents communicate locally to minimize the cost. The network topology is defined by a time-invariant doubly stochastic matrix $\Pb$, where $[\Pb]_{ji}>0$ if agent $j$ communicates with agent $i$; otherwise $[\Pb]_{ji}=0$. The network is assumed to be connected, and there exists a geometric mixing bound for $\Pb$ \cite{liu2008monte}, such that 
$$\sum_{j=1}^m \left|[\Pb^k]_{ji}-1/m\right|\leq \sqrt{m}\beta^k,\:i\in [m],$$
where $\beta$ is the second largest singular value of $\Pb$.

\subsection{Strong Stability and Sequential Strong Stability}\label{sec:stability}
We consider the set of strongly stable linear (i.e., $\ub=\Kb\xb$) controllers as the benchmark policy class. Following \cite{cohen2018online}, we define the notion of strong stability as follows. 
\begin{definition}\label{D: Strong Stability}(Strong Stability) A linear policy $\Kb$ is $(\kappa, \gamma)$-strongly stable (for $\kappa > 0$ and $0<\gamma\leq 1$) for the LTI system $(\Ab,\Bb)$, if $\norm{\Kb}\leq \kappa$, and there exist matrices $\Lb$ and $\Hb$ such that $\Ab+\Bb\Kb=\Hb\Lb\Hb^{-1}$, with $\norm{\Lb}\leq 1-\gamma$ and $\norm{\Hb}\norm{\Hb^{-1}}\leq \kappa$.
\end{definition}
Intuitively, a strongly stable policy ensures fast mixing
and exponential convergence to a steady-state distribution. In particular, for the LTI system $\xb_{t+1}=\Ab\xb_t + \Bb\ub_t + \wb_t$, if a $(\kappa, \gamma)$-strongly stable policy $\Kb$ is applied ($\ub_t = \Kb\xb_t$), $\widehat{\Xb}_t$ (the state covariance matrix of $\xb_t$) converges to $\Xb$ (the steady-state covariance matrix) with the following exponential rate $$\norm{\widehat{\Xb}_t-\Xb}\leq \kappa^2 e^{-2\gamma t}\norm{\widehat{\Xb}_0-\Xb}.$$
See Lemma 3.2 in \cite{cohen2018online} for details. The {\it sequential} nature of {\it online} LQ control requires another notion of strong stability, called {\it sequential strong stability} \citet{cohen2018online}, defined as follows.

\begin{definition}\label{D: Sequential Strong Stability}(Sequential Strong Stability) A sequence of linear  policies $\{\Kb_t\}_{t=1}^T$ is $(\kappa,\gamma)$-strongly stable if there exist matrices $\{\Hb_t\}_{t=1}^T$ and $\{\Lb_t\}_{t=1}^T$ such that $\Ab+\Bb\Kb_t=\Hb_t\Lb_t\Hb_t^{-1}$ for all $t$ with the following properties:
\begin{enumerate}
    \item $\norm{\Lb_t}\leq 1-\gamma$ and $\norm{\Kb_t}\leq \kappa$.
    \item $\norm{\Hb_t}\leq \beta$ and $\norm{\Hb_t^{-1}}\leq 1/\alpha$ with $\kappa=\beta/\alpha$.
    \item $\norm{\Hb_{t+1}^{-1}\Hb_t}\leq 1+\gamma/2$.
\end{enumerate}
\end{definition}
Sequential strong stability generalizes strong stability to the time-varying scenario. On the technical level, it helps with characterizing the convergence of the state covariance matrices when a sequence of policies $\{\Kb_t\}_{t=1}^T$ is used instead of a fixed policy $\Kb$, which will be the case in this work. 

\subsection{SDP Relaxation for LQ Control}
For the following dynamical system
$$\xb_{t+1} = \Ab\xb_{t} + \Bb\ub_{t} + \wb_{t},\quad \wb_{t}\sim \Nc(0,\Wb),$$
the infinite-horizon version of \eqref{eq:cost-central}, i.e., $$\minimize \lim_{T\rightarrow\infty}J_T(\pi)/T,$$
with fixed cost matrices $\Qb$ and $\Rb$ can be relaxed via a semi-definite programming when the steady-state distribution exists. For $\nu>0$, the SDP relaxation is formulated as \cite{cohen2018online}
\begin{equation}\label{eq:SDP}
\begin{split}
    \text{minimize}\quad &J(\Sigma)=\begin{pmatrix}\Qb & 0\\0 & \Rb\end{pmatrix}\bullet\Sigma\\
    \text{subject to}\quad &\Sigma_{\xb\xb}=(\Ab\:\Bb)\Sigma(\Ab\:\Bb)^\top+\Wb,\\
    &\Sigma\succeq 0,\quad \text{Tr}(\Sigma)\leq\nu,
\end{split}
\end{equation}
where $\Sigma=\begin{pmatrix}\Sigma_{\xb\xb} & \Sigma_{\xb\ub}\\\Sigma_{\ub\xb} & \Sigma_{\ub\ub}\end{pmatrix}$. Recall that in the online LQ problem, we deal with time-varying cost matrices $(\Qb_t,\Rb_t)$, and for any $t\in [T]$, the above SDP yields different solutions. 

In fact, for any feasible solution $\Sigma$ of the above SDP, a strongly stable controller $\Kb=\Sigma_{\xb\ub}^\top \Sigma_{\xb\xb}^{-1}$ can be extracted. The steady-state covariance matrix induced by this controller is also feasible for the SDP and its cost is at most that of $\Sigma$ (see Theorem 4.2 in \cite{cohen2018online}).

Moreover, for any (slowly-varying) sequence of feasible solutions to the above SDP, the induced controller sequence is sequentially strongly-stable. This implies that the covariance matrix of the state converges to the steady-state in a rapid sense as the following.

\begin{lemma}\label{L: Convergence of Sequential Strong Stability (feasible solutions)} (Lemma 4.4 in \cite{cohen2018online}) Assume that $\Wb\succeq \sigma^2\Ib$ and let $\kappa=\sqrt{\nu}/\sigma$. Let $\{\Sigma_t\}$ be a sequence of feasible solutions of the SDP, and suppose that $\norm{\Sigma_{t+1}-\Sigma_t}\leq\eta$ for all $t$ and for some $\eta\leq \sigma^2/\kappa^2$. Then, the control matrix $\Kb_t=(\Sigma_{t})_{\xb\ub}^\top(\Sigma_{t})_{\xb\xb}^{-1}$ is $(\kappa,\frac{1}{2\kappa^2})$-strongly stable for all $t$. 
\end{lemma}

Furthermore, it can be shown that given the sequence $\Xb_t=(\Sigma_{t})_{\xb\xb}$, if we follow the policy sequence $\pi_t(\xb)=\Kb_t\xb+v_t$ where $v_t\sim \Nc\left(0,(\Sigma_{t})_{\ub \ub}-\Kb_{t}(\Sigma_{t})_{\xb \xb}\Kb_{t}^\top\right)$, the following relationship holds:
$$\norm{\widehat{\Xb}_{t+1}-\Xb_{t+1}}\leq \kappa^2 e^{-(\frac{1}{2\kappa^2}) t}\norm{\widehat{\Xb}_1-\Xb_1}+4\eta\kappa^4,$$
where $\widehat{\Xb}_t$ is the state covariance matrix on round $t$ \cite{cohen2018online}.

\section{Algorithm and Theoretical Results}
We now lay out the distributed online LQ algorithm and provide its theoretical regret bound. 
\subsection{Algorithm}
In the distributed online LQ, each agent $i$ at time $t$ maintains an ideal steady-state covariance matrix $\Sigma_{i,t}$ by running a distributed online gradient descent on the SDP \eqref{eq:SDP}. Then, a control matrix $\Kb_{i,t}$ is extracted from $\Sigma_{i,t}$ and is used to determine the action. In particular, the action $\ub_{i,t}$ is sampled from a Gaussian distribution $\Nc(\Kb_{i,t}\xb_{i,t}, \Vb_{i,t})$, which entails $\mathrm{E}[\ub_{i,t}|\Fc_t]=\Kb_{i,t}\xb_{i,t}$, where $\Fc_t$ is the smallest $\sigma$-field containing the
information about all agents up to time $t$. This stochastic policy ensures the fast convergence of the covariance matrix of $\xb_{i,t}$ and $\ub_{i,t}$ to the iterate $\Sigma_{i,t}$ generated by the algorithm. 
The proposed method is outlined in Algorithm \ref{alg:Online Distributed LQ Control}.

\subsection{Theoretical Result: Regret Bound}
In this section, we present our main theoretical result. By applying algorithm \ref{alg:Online Distributed LQ Control}, we show that for a multi-agent network of identical LTI systems (with a connected communication graph), the individual regret of an arbitrary agent is upper-bounded by $O(\sqrt{T})$, which implies that the performance of all agents would converge to that of the best fixed controller in hindsight for large enough $T$.\\
\begin{algorithm}[tb]
   \caption{Online Distributed LQ Control}
   \label{alg:Online Distributed LQ Control}
\begin{algorithmic}[1]
   \STATE {\bfseries Require:} number of agents $m$, doubly stochastic matrix $\Pb\in \mathrm{R}^{m\times m}$, parameter $\nu$, step size $\eta$, system matrices $(\Ab,\Bb)$.
   
  \STATE {\bf Initialize:} $\Sigma_{i,1}$ is identically initialized with a feasible point and $\xb_{i,1}$ is drawn from normal distribution with mean zero for $i\in [m]$.
  
   \FOR{$t=1,2,\ldots,T$}
        \FOR{$i=1,2,\ldots,m$}
        \STATE Receive $\xb_{i,t}$
        \STATE Compute $\Kb_{i,t}=(\Sigma_{i,t})_{\ub \xb}(\Sigma_{i,t})_{\xb \xb}^{-1}, \Vb_{i,t}=(\Sigma_{i,t})_{\ub \ub}-\Kb_{i,t}(\Sigma_{i,t})_{\xb \xb}\Kb_{i,t}^\top$
        \STATE Predict $\ub_{i,t}\sim \Nc(\Kb_{i,t}\xb_{i,t}, \Vb_{i,t})$ and Observe $\Qb_{i,t},\Rb_{i,t}$
        \STATE Communicate $\Sigma_{i,t}$ with agents in the neighborhood and obtain their parameters
        \STATE $\Sigma_{i,t+1} = \Pi_{\Sc}\left[\sum_{j}\Pb_{ji}\Sigma_{j,t}-\eta \begin{pmatrix}\Qb_{i,t} & 0\\0 & \Rb_{i,t}
        \end{pmatrix}\right]$,
        where 
            \begin{equation*}
            \resizebox{.85\hsize}{!}{
            $\Sc=\left\{\Sigma\in\mathrm{R}^{
            (d+k)\times (d+k)}\middle| 
            \begin{tabular}{@{}l@{}}$\Sigma \succeq 0,\quad \text{Tr}(\Sigma)\leq \nu,$\\ $\Sigma_{\xb\xb}=(\Ab\: \Bb)\Sigma(\Ab\:\Bb)^\top+\Wb$\end{tabular}
            \right\}$} 
            \end{equation*}
            
        \ENDFOR
   \ENDFOR
\end{algorithmic}
\end{algorithm}

\begin{theorem}\label{T: Online Distributed LQ Controller}
Assume that the network is connected, $\text{Tr}(\Wb)\leq \lambda^2$ and $\Wb\succeq\sigma^2\Ib$. Given $\kappa>1$ and $0\leq\gamma<1$, set $\nu=2\kappa^4\lambda^2/\gamma$ and $\eta=1/\sqrt{\rho T}$, where $$\rho=\left[4mC^2\left(3+\frac{4\sqrt{m}}{1-\beta}\right) + mC(1+\frac{\nu}{\sigma^2})\frac{16\sqrt{2m}C\nu}{(1-\beta)\sigma^2}\right].$$
By running Algorithm \ref{alg:Online Distributed LQ Control}, the expected individual regret of agent $j$ with respect to any $(\kappa,\gamma)$-strongly stable controller $\Kb^s$ is bounded as follows
$$\text{Regret}_T^j(\mathcal{A})=J^j_T(\mathcal{A})-J_T(\Kb^s)=O\left((1-\beta)^{-0.5}\sqrt{T}\right),$$
for $T\geq\left(\frac{4\sqrt{2}\nu C}{\sigma^4(1-\beta)\rho^{1/2}}\right)^2$.
\end{theorem}

The dependence of regret bound to the spectral gap $1-\beta$ is perhaps not surprising, as it has been previously observed in distributed online algorithms (see e.g., \cite{shahrampour2018distributed} Corollary 4).


\section{Numerical Experiments}
\begin{figure}[t] 
    \includegraphics[width=1\columnwidth]{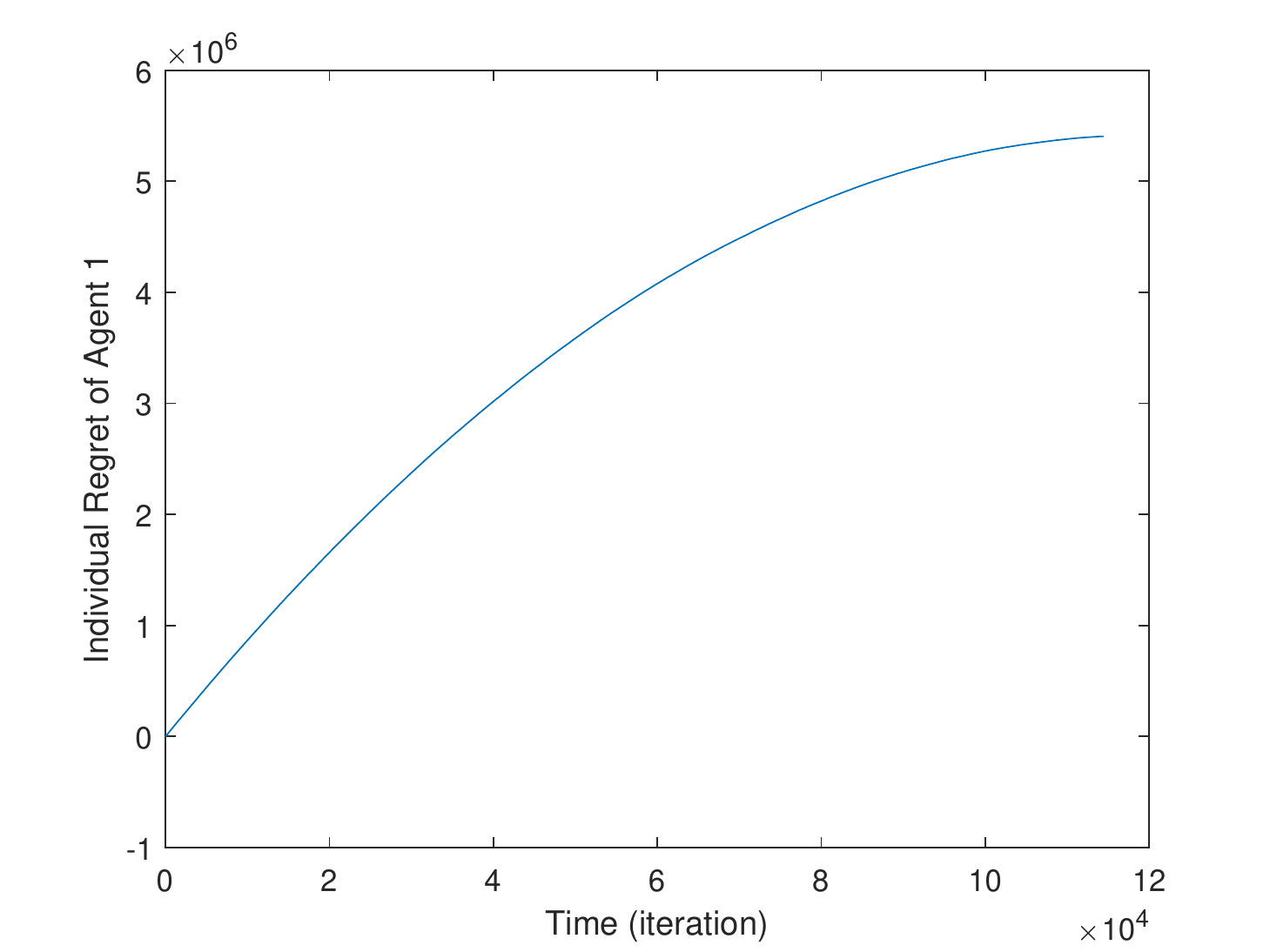}
    \caption{The plot of individual regret of agent 1 vs. time shows the sub-linearity of the regret evolvement.}
    \label{fig:individual regret}
\end{figure}

\begin{figure}[t!] 
    \includegraphics[width=1\columnwidth]{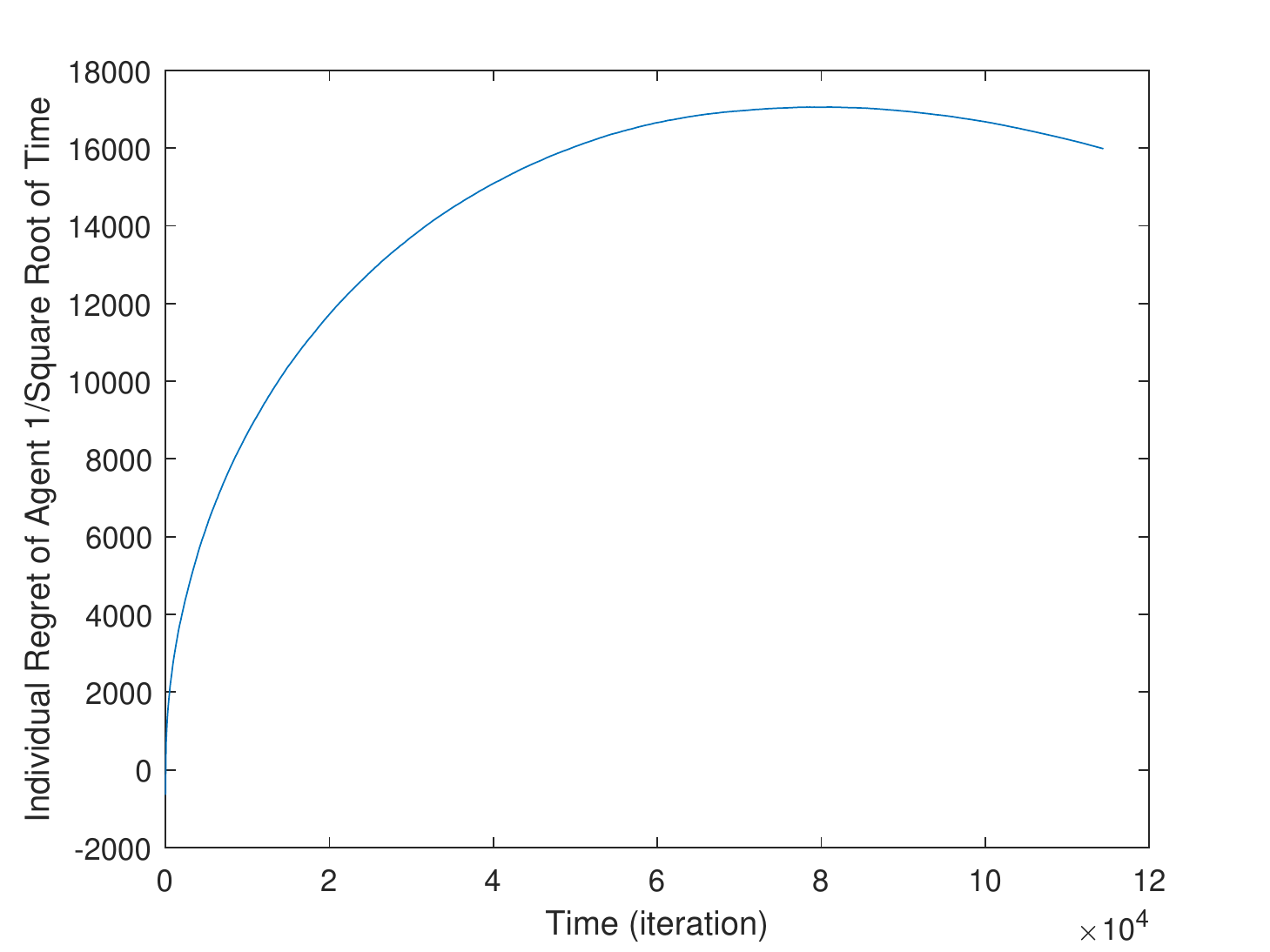}
    \caption{This plot shows that the individual regret of agent 1 is of $O(\sqrt{T})$ when $T$ is large enough.}
    \label{fig:individual regret sqrt}
\end{figure}
We now provide numerical simulations verifying the theoretical guarantee of our algorithm.

\vspace{0.2cm}
\noindent
\textbf{Experiment Setup:} We consider a distributed network of five agents where $d=k=3$. The network topology is a cycle, where each agent has a self-weight of 0.6, and the rest of the weight is evenly distributed between its neighborhood. The resulting communication matrix is as follows
$$\Pb=\begin{pmatrix}0.6 & 0.2 & 0 & 0 & 0.2\\
                 0.2 & 0.6 & 0.2 & 0 & 0\\
                 0 & 0.2 & 0.6 & 0.2 & 0\\
                 0 & 0 & 0.2 & 0.6 & 0.2\\
                 0.2 & 0 & 0 & 0.2 & 0.6
\end{pmatrix}.$$
The other hyper-parameters are set as follows: $\kappa=1.5$, $\gamma=0.4$, $C=30$. We let matrices $\Ab=(1-2\gamma)\Ib$ and $\Bb=(\gamma/\kappa)\Ib$. We set the cost matrix $\Qb_{i,t}$ (respectively, $\Rb_{i,t}$) as a diagonal matrix with each diagonal entry sampled from the uniform distribution over $[0,C/d]$ (respectively, $[0,C/k]$) to ensure that $\text{Tr}(\Qb_{i,t}),\text{Tr}(\Rb_{i,t})\leq C$. The noise $\wb_{i,t}$ is sampled from a standard Gaussian distribution, and thus $\lambda^2=d=3$ and $\sigma^2=1$.

\vspace{0.2cm}
\noindent
\textbf{Simulation:} The total iteration number $T$ is set as 30 times of the theoretical lower bound in Theorem \ref{T: Online Distributed LQ Controller} in order to better see the performance. We let $\Kb^s=(1e-2)(-\kappa)\Ib$ which is $(\kappa,\gamma)$-strongly stable with $\Ab,\Bb$, and leads to a small enough cumulative cost to be the benchmark. Noting that we apply Dykstra's projection algorithm for the projection step, the matrix $\Vb_{i,t}$ for action-sampling may not be positive semi-definite (PSD) due to floating-point computations, so we do some tuning by adding to it a small term ($(1e-25)\Ib$) to keep it PSD. The parameters $\Sigma_{i,1}$ are identically initialized and the initial states of all agents are sampled from normal distribution. The entire process is repeated for 30 Monte-Carlo simulations. 
\begin{figure}[t!] 
    \includegraphics[width=1\columnwidth]{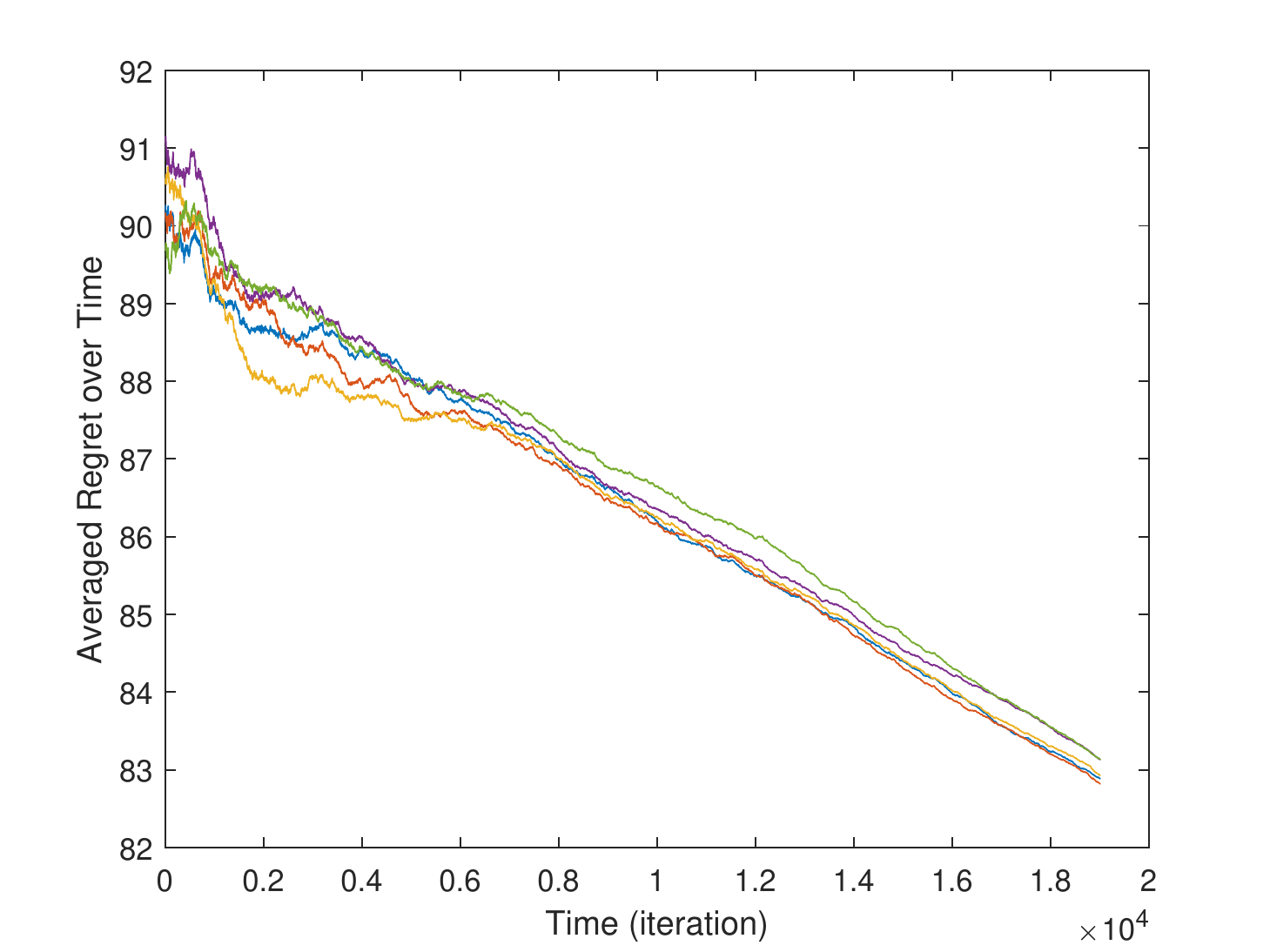}
    \caption{The averaged regrets over time for all agents converge as time grows.}
    \label{fig:consensus}
\end{figure}

\vspace{0.2cm}
\noindent
\textbf{Performance:} To see the sub-linearity of the individual regret established in Theorem \ref{T: Online Distributed LQ Controller}, we plot the regret of agent 1 over time. Fig. \ref{fig:individual regret} shows that the regret is upper-bounded by a linear function. To better capture this, we also plot the regret normalized by the root-square of time in Fig. \ref{fig:individual regret sqrt}. we observe that for large enough $T$, the slope of the plot is non-positive, which verifies that the regret is upper-bounded by $O(\sqrt{T})$. We also present the performance of all agents to verify that they reach consensus asymptotically. In Fig. \ref{fig:consensus}, it can be seen that the time-averaged regrets for all agents are decreasing and converging together. 

\section{Conclusion}
In this paper, we considered the distributed online LQ problem with known identical LTI systems and decoupled, time-varying quadratic cost functions. We developed a fully distributed algorithm to minimize the finite-horizon cost, which can be recast as a regret minimization. We proved that the individual regret, which is the performance of the control sequence of any agent compared to the best (linear and strongly stable) controller in hindsight, is upper bounded by $O(\sqrt{T})$. Our numerical simulations verified that the theoretical bound is indeed valid. For future works, it is important to analyze the distributed online LQ under more general settings. Possible directions include extending the setup to {\it unknown} dynamics or assuming {\it coupled} time-varying cost functions. 

\appendix
{{\bf\itshape Proof of Theorem \ref{T: Online Distributed LQ Controller}: }} Recalling the definition of regret \eqref{eq:regret}, for a fixed arbitrary $(\kappa,\gamma)$-strongly stable controller $\Kb^s$ and agent $j$, the regret is expressed as the following:
\begin{equation}\label{eq: Regret Proof eq1}
\begin{split}
    &J^j_T(\mathcal{A})-J_T(\Kb^s)\\
    =&\mathrm{E}\left[\sum_{t=1}^T\sum_{i=1}^m(\xb_{j,t}^\top\Qb_{i,t}\xb_{j,t} + \ub_{j,t}^\top\Rb_{i,t}\ub_{j,t})\right]\\
    -&\mathrm{E}\left[\sum_{t=1}^T(\xb_{t}^{s\top}\Qb_{t}\xb_{t}^s + \ub_{t}^{s\top}\Rb_{t}\ub_{t}^s)\right],
\end{split}    
\end{equation}
where $\ub_t^s = \Kb^s\xb_t^s$ for all $t$.

Let us denote 
$$
\Lb_{i,t}=\begin{pmatrix}\Qb_{i,t}&0\\0&\Rb_{i,t}\end{pmatrix}~~~~\text{and}~~~~\Lb_{t}=\begin{pmatrix}\Qb_{t}&0\\0&\Rb_{t}\end{pmatrix},
$$
where $\Lb_{t}=\sum_{i=1}^m\Lb_{i,t}$. Also let, 
\begin{align*}
    \widehat{\Sigma}_{j,t}&=\mathrm{E}\left[[\xb_{j,t}^\top\:\ub_{j,t}^\top]^\top[\xb_{j,t}^\top\:\ub_{j,t}^\top]\right]\\
    \widehat{\Sigma}_t^s&=\mathrm{E}\left[[\xb_t^{s\top}\:\ub_t^{s\top}]^\top[\xb_t^{s\top}\:\ub_t^{s\top}]\right].
\end{align*}
We can then write \eqref{eq: Regret Proof eq1} as
\begin{equation}\label{eq: Regret Proof eq2}
\begin{split}
    &\sum_{t=1}^T\sum_{i=1}^m\Lb_{i,t}\bullet\widehat{\Sigma}_{j,t} - \sum_{t=1}^T\Lb_t\bullet\widehat{\Sigma}_t^s\\
    =&\sum_{t=1}^T\sum_{i=1}^m\Lb_{i,t}\bullet(\widehat{\Sigma}_{j,t}-\Sigma_{j,t})\\ + &\sum_{t=1}^T\sum_{i=1}^m\Lb_{i,t}\bullet\Sigma_{j,t} - \sum_{t=1}^T\Lb_t\bullet\Sigma^s\\ + &\sum_{t=1}^T\Lb_t\bullet(\Sigma^s-\widehat{\Sigma}_t^s),
\end{split}    
\end{equation}
where $\Sigma^s$ is the steady-state covariance matrix induced by $\Kb^s$, and $\Sigma_{j,t}$ is generated by Algorithm \ref{alg:Online Distributed LQ Control}. Now, we show how each term in \eqref{eq: Regret Proof eq2} is bounded.\\\\
{\bf (I)} For the term $\sum_{t=1}^T\sum_{i=1}^m\Lb_{i,t}\bullet\Sigma_{j,t} - \sum_{t=1}^T\Lb_t\bullet\Sigma^s$:\\
Based on Lemma 3.3 in \cite{cohen2018online}, it can be shown that $\text{Tr}(\Sigma^s)=\text{Tr}(\Sigma^s_{\xb\xb})+\text{Tr}(\Sigma^s_{\ub\ub})\leq 2\kappa^4\lambda^2/\gamma=\nu$. Then, by Lemma 4.1 in \cite{cohen2018online}, $\Sigma^s$ is a feasible solution to the SDP \eqref{eq:SDP}. Based on the definition of the feasible set $\Sc$, the diameter $\sup_{\Sigma,\Sigma^\prime\in\Sc}\norm{\Sigma-\Sigma^\prime}_F\leq 2\sup_{\Sigma\in\Sc}\norm{\Sigma}_F=2\sup_{\Sigma\in\Sc}\sqrt{\text{Tr}(\Sigma^2)}\leq 2\sup_{\Sigma\in\Sc}\sqrt{\text{Tr}(\Sigma)^2}\leq 2\nu$. And the norm of the gradient of the linear loss function $\Sc\to\Lb_{i,t}\bullet\Sc$ is upper bounded by $\sqrt{2}C$ since $\sqrt{\text{Tr}(\Qb_{i,t}^\top\Qb_{i,t}) + \text{Tr}(\Rb_{i,t}^\top\Rb_{i,t})}\leq \sqrt{2}C$.

Let $\Sigma^*=\argmin_{\Sigma\in\Sc}\sum_{t=1}^T\Lb_t\bullet\Sigma$. Based on the regret bound of distributed online gradient descent \cite{yan2012distributed}, we have
\begin{equation}\label{eq: Regret Proof eq3}
\begin{split}
    &\sum_{t=1}^T\sum_{i=1}^m\Lb_{i,t}\bullet\Sigma_{j,t} - \sum_{t=1}^T\Lb_t\bullet\Sigma^s\\
    \leq &\sum_{t=1}^T\sum_{i=1}^m\Lb_{i,t}\bullet\Sigma_{j,t} - \sum_{t=1}^T\Lb_t\bullet\Sigma^*\\
    \leq &\frac{m\nu}{\eta} + \left(3+\frac{4\sqrt{m}}{1-\beta}\right)4mC^2\eta T,
\end{split}    
\end{equation}
where $\beta \in [0,1)$ is the second largest singular value of $\Pb$. Also, based on Lemma 3 in \cite{yan2012distributed}, the variation $\norm{\Sigma_{j,t+1}-\Sigma_{j,t}}_F$ is upper bounded as the following:
\begin{equation}\label{eq: Regret Proof eq4}
    \norm{\Sigma_{j,t+1}-\Sigma_{j,t}}_F\leq \frac{4\sqrt{2m}C\eta}{1-\beta}.
\end{equation}
{\bf (II)} For the term $\sum_{t=1}^T\sum_{i=1}^m\Lb_{i,t}\bullet(\widehat{\Sigma}_{j,t}-\Sigma_{j,t})$:\\
Based on Algorithm \ref{alg:Online Distributed LQ Control}, we have 
\begin{equation*}
\begin{split}
    \Sigma_{j,t} = & \begin{pmatrix}(\Sigma_{j,t})_{\xb\xb}&(\Sigma_{j,t})_{\xb\ub}\\(\Sigma_{j,t})_{\ub\xb}&(\Sigma_{j,t})_{\ub\ub}\end{pmatrix}\\
    =&\begin{pmatrix}(\Sigma_{j,t})_{\xb\xb}&(\Sigma_{j,t})_{\xb\xb}\Kb_{j,t}^\top\\\Kb_{j,t}(\Sigma_{j,t})_{\xb\xb}&\Kb_{j,t}(\Sigma_{j,t})_{\xb\xb}\Kb_{j,t}^\top\end{pmatrix}\\ 
    +&\begin{pmatrix}0&0\\0&\Vb_{j,t}\end{pmatrix} 
\end{split}
\end{equation*}
and
\begin{equation*}
\begin{split}
    &\widehat{\Sigma}_{j,t}\\
    =&\begin{pmatrix}(\widehat{\Sigma}_{j,t})_{\xb\xb}&(\widehat{\Sigma}_{j,t})_{\xb\xb}\Kb_{j,t}^\top\\\Kb_{j,t}(\widehat{\Sigma}_{j,t})_{\xb\xb}&\Kb_{j,t}(\widehat{\Sigma}_{j,t})_{\xb\xb}\Kb_{j,t}^\top\end{pmatrix} + \begin{pmatrix}0&0\\0&\Vb_{j,t}\end{pmatrix}.
\end{split}    
\end{equation*}
Therefore, we get
\begin{equation}\label{eq: Regret Proof eq5}
\begin{split}
    &\Lb_{i,t}\bullet(\widehat{\Sigma}_{j,t}-\Sigma_{j,t})\\
    =&\Qb_{i,t}\bullet\big((\widehat{\Sigma}_{j,t})_{\xb\xb}-(\Sigma_{j,t})_{\xb\xb}\big)\\
    +&\Rb_{i,t}\bullet\Kb_{j,t}\big((\widehat{\Sigma}_{j,t})_{\xb\xb}-(\Sigma_{j,t})_{\xb\xb}\big)\Kb_{j,t}^\top\\
    =&(\Qb_{i,t} + \Kb_{j,t}^\top\Rb_{i,t}\Kb_{j,t})\bullet\big((\widehat{\Sigma}_{j,t})_{\xb\xb}-(\Sigma_{j,t})_{\xb\xb}\big)\\
    \leq &\text{Tr}(\Qb_{i,t} + \Kb_{j,t}^\top\Rb_{i,t}\Kb_{j,t})\norm{(\widehat{\Sigma}_{j,t})_{\xb\xb}-(\Sigma_{j,t})_{\xb\xb}}\\
    \leq &\left[\text{Tr}(\Qb_{i,t}) + \text{Tr}(\Rb_{i,t})\norm{\Kb_{j,t}\Kb_{j,t}^\top}\right]\norm{(\widehat{\Sigma}_{j,t})_{\xb\xb}-(\Sigma_{j,t})_{\xb\xb}}\\
    \leq &C(1+\frac{\nu}{\sigma^2})\norm{(\widehat{\Sigma}_{j,t})_{\xb\xb}-(\Sigma_{j,t})_{\xb\xb}},
\end{split}    
\end{equation}
where the third inequality holds since $\text{Tr}(\Qb_{i,t}),\text{Tr}(\Rb_{i,t})\leq C$ and $\Kb_{j,t}$ is $(\frac{\sqrt{\nu}}{\sigma},\frac{\sigma^2}{2\nu})$-strongly stable based on Lemma 4.3 in \cite{cohen2018online}.\\\\
Choosing $\eta$ such that $\frac{4\sqrt{2m}C\eta}{1-\beta}\leq \frac{\sigma^4}{\nu}$, based on \eqref{eq: Regret Proof eq4} and Lemma \ref{L: Convergence of Sequential Strong Stability (feasible solutions)}, we have
\begin{equation}\label{eq: Regret Proof eq6}
\begin{split}
    &\norm{(\widehat{\Sigma}_{j,t})_{\xb\xb} - (\Sigma_{j,t})_{\xb\xb}}\\
    \leq &\frac{\nu}{\sigma^2}e^{-\frac{\sigma^2}{2\nu}(t-1)}\norm{(\widehat{\Sigma}_{j,1})_{\xb\xb} - (\Sigma_{j,1})_{\xb\xb}} + \frac{16\sqrt{2m}C\eta\nu}{(1-\beta)\sigma^2}.
\end{split}
\end{equation}
Substituting \eqref{eq: Regret Proof eq6} into \eqref{eq: Regret Proof eq5} and summing over $t\in[T]$, we get
{\small
\begin{equation}\label{eq: Regret Proof eq7}
\begin{split}
    &\sum_{t=1}^T\Lb_{i,t}\bullet(\widehat{\Sigma}_{j,t}-\Sigma_{j,t})\\
    \leq &C(1+\frac{\nu}{\sigma^2})\left(\frac{\nu}{\sigma^2}\norm{(\widehat{\Sigma}_{j,1})_{\xb\xb} - (\Sigma_{j,1})_{\xb\xb}}\sum_{t=1}^Te^{-\frac{\sigma^2}{2\nu}(t-1)}\right)\\
    + &C(1+\frac{\nu}{\sigma^2})\left(\frac{16\sqrt{2m}C\eta\nu}{(1-\beta)\sigma^2}T\right)\\
    \leq &C(1+\frac{\nu}{\sigma^2})\left(\frac{2\nu^2}{\sigma^4}\norm{(\widehat{\Sigma}_{j,1})_{\xb\xb} - (\Sigma_{j,1})_{\xb\xb}} + \frac{16\sqrt{2m}C\eta\nu}{(1-\beta)\sigma^2}T\right),
\end{split}
\end{equation}}
where the second inequality comes from the fact that $\sum_{t=1}^Te^{-\alpha t}\leq \int_0^\infty e^{-\alpha t}dt=1/\alpha$ for $\alpha>0$. Summing \eqref{eq: Regret Proof eq7} over $i$, the result is obtained.\\\\
{\bf (III)} For the term $\sum_{t=1}^T\Lb_t\bullet(\Sigma^s-\widehat{\Sigma}_t^s)$:\\
By denoting $\Sigma^s = \begin{pmatrix}\Sigma^s_{\xb\xb} & \Sigma^s_{\xb\xb}\Kb^{s\top}\\\Kb^s\Sigma^s_{\xb\xb}&\Kb^s\Sigma^s_{\xb\xb}\Kb^{s\top}\end{pmatrix}$ and $\widehat{\Sigma}^s_t = \begin{pmatrix}(\widehat{\Sigma}^s_t)_{\xb\xb} & (\widehat{\Sigma}^s_t)_{\xb\xb}\Kb^{s\top}\\\Kb^s(\widehat{\Sigma}^s_t)_{\xb\xb}&\Kb^s(\widehat{\Sigma}^s_t)_{\xb\xb}\Kb^{s\top}\end{pmatrix}$, we have
\begin{equation}\label{eq: Regret Proof eq8}
\begin{split}
    &\Lb_t\bullet(\Sigma^s-\widehat{\Sigma}^s_t)\\
    =&\sum_{i=1}^m(\Qb_{i,t}+\Kb^s\Rb_{i,t}\Kb^{s\top})\bullet\left(\Sigma^s_{\xb\xb}-(\widehat{\Sigma}^s_t)_{\xb\xb}\right)\\
    \leq &\sum_{i=1}^m \text{Tr}(\Qb_{i,t}+\Kb^s\Rb_{i,t}\Kb^{s\top})\norm{\Sigma^s_{\xb\xb}-(\widehat{\Sigma}^s_t)_{\xb\xb}}\\
    \leq &mC(1+\kappa^2)\norm{\Sigma^s_{\xb\xb}-(\widehat{\Sigma}^s_t)_{\xb\xb}},
\end{split}
\end{equation}
where the second inequality comes from the fact that $\text{Tr}(\Qb_{i,t}),\text{Tr}(\Rb_{i,t})\leq C$ and $\Kb^s$ is $(\kappa,\gamma)$-strongly stable.\\\\
Based on Lemma 3.2 in \cite{cohen2018online}, we get
\begin{equation}\label{eq: Regret Proof eq9}
    \norm{(\widehat{\Sigma}^s_t)_{\xb\xb}-\Sigma^s_{\xb\xb}}\leq \kappa^2 e^{-2\gamma(t-1)}\norm{(\widehat{\Sigma}^s_1)_{\xb\xb}-\Sigma^s_{\xb\xb}}.
\end{equation}
Substituting \eqref{eq: Regret Proof eq9} into \eqref{eq: Regret Proof eq8} and summing over $t$, we have 
\begin{equation}\label{eq: Regret Proof eq10}
\begin{split}
    &\sum_{t=1}^T\Lb_t\bullet(\Sigma^s - \widehat{\Sigma}_t^s)\\
    \leq &mC(1+\kappa^2)\kappa^2\norm{(\widehat{\Sigma}^s_1)_{\xb\xb}-\Sigma^s_{\xb\xb}}\sum_{t=1}^Te^{-2\gamma(t-1)}\\
    \leq &\frac{mC(\kappa^2 + \kappa^4)}{2\gamma}\norm{(\widehat{\Sigma}^s_1)_{\xb\xb}-\Sigma^s_{\xb\xb}}.
\end{split}
\end{equation}
Based on \eqref{eq: Regret Proof eq3}, \eqref{eq: Regret Proof eq7} and \eqref{eq: Regret Proof eq10}, we have 
\begin{equation}\label{eq: Regret Proof eq11}
\begin{split}
    &J^j_T(\mathcal{A})-J_T(\Kb^s)\\
    \leq &\frac{m\nu}{\eta} + mC(\frac{2\nu^2(\sigma^2+\nu)}{\sigma^6})\norm{(\widehat{\Sigma}_{j,1})_{\xb\xb}-(\Sigma_{j,1})_{\xb\xb}}\\
    + &\frac{mC(\kappa^2 + \kappa^4)}{2\gamma}\norm{(\widehat{\Sigma}^s_1)_{\xb\xb}-\Sigma^s_{\xb\xb}} + \rho\eta T,
\end{split}
\end{equation}
where $$\rho\triangleq\left[4mC^2\left(3+\frac{4\sqrt{m}}{1-\beta}\right) + mC(1+\frac{\nu}{\sigma^2})\frac{16\sqrt{2m}C\nu}{(1-\beta)\sigma^2}\right].$$

By setting $\eta=1/\sqrt{\rho T}$, the upper bound in \eqref{eq: Regret Proof eq11} is $O(\sqrt{\rho T})$. (Here it is assumed that $T\geq \left(\frac{4\sqrt{2m}\nu C}{\sigma^4(1-\beta)\rho^{1/2}}\right)^2$ to make sure $\frac{4\sqrt{2m}C\eta}{1-\beta}\leq \frac{\sigma^4}{\nu}$ and \eqref{eq: Regret Proof eq6} holds.)
\endproof
\addtolength{\textheight}{-0cm}



\bibliographystyle{IEEEtran}
\bibliography{references}

\begin{thebibliography}{10}
\providecommand{\url}[1]{#1}
\csname url@samestyle\endcsname
\providecommand{\newblock}{\relax}
\providecommand{\bibinfo}[2]{#2}
\providecommand{\BIBentrySTDinterwordspacing}{\spaceskip=0pt\relax}
\providecommand{\BIBentryALTinterwordstretchfactor}{4}
\providecommand{\BIBentryALTinterwordspacing}{\spaceskip=\fontdimen2\font plus
\BIBentryALTinterwordstretchfactor\fontdimen3\font minus
  \fontdimen4\font\relax}
\providecommand{\BIBforeignlanguage}[2]{{%
\expandafter\ifx\csname l@#1\endcsname\relax
\typeout{** WARNING: IEEEtran.bst: No hyphenation pattern has been}%
\typeout{** loaded for the language `#1'. Using the pattern for}%
\typeout{** the default language instead.}%
\else
\language=\csname l@#1\endcsname
\fi
#2}}
\providecommand{\BIBdecl}{\relax}
\BIBdecl

\bibitem{4309169}
B.~D.~O. {Anderson}, J.~B. {Moore}, and B.~P. {Molinari}, ``Linear optimal
  control,'' \emph{IEEE Transactions on Systems, Man, and Cybernetics}, vol.
  SMC-2, no.~4, pp. 559--559, 1972.

\bibitem{bertsekas1995dynamic}
D.~P. Bertsekas, \emph{Dynamic programming and optimal control}, vol.~1, no.~2.

\bibitem{zhou1996robust}
K.~Zhou, J.~C. Doyle, K.~Glover \emph{et~al.}, \emph{Robust and optimal
  control}.\hskip 1em plus 0.5em minus 0.4em\relax Prentice hall New Jersey,
  1996, vol.~40.

\bibitem{cohen2018online}
A.~Cohen, A.~Hasidim, T.~Koren, N.~Lazic, Y.~Mansour, and K.~Talwar, ``Online
  linear quadratic control,'' in \emph{International Conference on Machine
  Learning}, 2018, pp. 1029--1038.

\bibitem{yan2012distributed}
F.~Yan, S.~Sundaram, S.~Vishwanathan, and Y.~Qi, ``Distributed autonomous
  online learning: Regrets and intrinsic privacy-preserving properties,''
  \emph{IEEE Transactions on Knowledge and Data Engineering}, vol.~25, no.~11,
  pp. 2483--2493, 2012.

\bibitem{4626964}
F.~{Borrelli} and T.~{Keviczky}, ``Distributed lqr design for identical
  dynamically decoupled systems,'' \emph{IEEE Transactions on Automatic
  Control}, vol.~53, no.~8, pp. 1901--1912, 2008.

\bibitem{6862471}
A.~{Mosebach} and J.~{Lunze}, ``Synchronization of autonomous agents by an
  optimal networked controller,'' in \emph{2014 European Control Conference
  (ECC)}, 2014, pp. 208--213.

\bibitem{5299181}
Y.~{Cao} and W.~{Ren}, ``Optimal linear-consensus algorithms: An lqr
  perspective,'' \emph{IEEE Transactions on Systems, Man, and Cybernetics, Part
  B (Cybernetics)}, vol.~40, no.~3, pp. 819--830, 2010.

\bibitem{8736845}
J.~{Jiao}, H.~L. {Trentelman}, and M.~K. {Camlibel}, ``A suboptimality approach
  to distributed linear quadratic optimal control,'' \emph{IEEE Transactions on
  Automatic Control}, vol.~65, no.~3, pp. 1218--1225, 2020.

\bibitem{8734804}
------, ``Distributed linear quadratic optimal control: Compute locally and act
  globally,'' \emph{IEEE Control Systems Letters}, vol.~4, no.~1, pp. 67--72,
  2020.

\bibitem{alemzadeh2019distributed}
S.~Alemzadeh and M.~Mesbahi, ``Distributed q-learning for dynamically decoupled
  systems,'' in \emph{2019 American Control Conference (ACC)}.\hskip 1em plus
  0.5em minus 0.4em\relax IEEE, 2019, pp. 772--777.

\bibitem{fattahi2019efficient}
S.~Fattahi, N.~Matni, and S.~Sojoudi, ``Efficient learning of distributed
  linear-quadratic controllers,'' \emph{arXiv preprint arXiv:1909.09895}, 2019.

\bibitem{furieri2020learning}
L.~Furieri, Y.~Zheng, and M.~Kamgarpour, ``Learning the globally optimal
  distributed lq regulator,'' in \emph{Learning for Dynamics and Control},
  2020, pp. 287--297.

\bibitem{furieri2020first}
L.~Furieri and M.~Kamgarpour, ``First order methods for globally optimal
  distributed controllers beyond quadratic invariance,'' in \emph{2020 American
  Control Conference (ACC)}.\hskip 1em plus 0.5em minus 0.4em\relax IEEE, 2020,
  pp. 4588--4593.

\bibitem{fazel2018global}
M.~Fazel, R.~Ge, S.~Kakade, and M.~Mesbahi, ``Global convergence of policy
  gradient methods for the linear quadratic regulator,'' in \emph{International
  Conference on Machine Learning}, 2018, pp. 1467--1476.

\bibitem{malik2019derivative}
D.~Malik, A.~Pananjady, K.~Bhatia, K.~Khamaru, P.~Bartlett, and M.~Wainwright,
  ``Derivative-free methods for policy optimization: Guarantees for linear
  quadratic systems,'' in \emph{The 22nd International Conference on Artificial
  Intelligence and Statistics}.\hskip 1em plus 0.5em minus 0.4em\relax PMLR,
  2019, pp. 2916--2925.

\bibitem{9130755}
H.~{Mohammadi}, M.~{Soltanolkotabi}, and M.~R. {Jovanović}, ``On the linear
  convergence of random search for discrete-time lqr,'' \emph{IEEE Control
  Systems Letters}, vol.~5, no.~3, pp. 989--994, 2021.

\bibitem{9147749}
H.~{Mohammadi}, M.~{Soltanolkotabi}, and M.~R. {Jovanovic}, ``Random search for
  learning the linear quadratic regulator,'' in \emph{2020 American Control
  Conference (ACC)}, 2020, pp. 4798--4803.

\bibitem{dean2019sample}
S.~Dean, H.~Mania, N.~Matni, B.~Recht, and S.~Tu, ``On the sample complexity of
  the linear quadratic regulator,'' \emph{Foundations of Computational
  Mathematics}, pp. 1--47, 2019.

\bibitem{hazan2017learning}
E.~Hazan, K.~Singh, and C.~Zhang, ``Learning linear dynamical systems via
  spectral filtering,'' in \emph{Advances in Neural Information Processing
  Systems}, 2017, pp. 6702--6712.

\bibitem{arora2018towards}
S.~Arora, E.~Hazan, H.~Lee, K.~Singh, C.~Zhang, and Y.~Zhang, ``Towards
  provable control for unknown linear dynamical systems,'' 2018.

\bibitem{agarwal2019online}
N.~Agarwal, B.~Bullins, E.~Hazan, S.~M. Kakade, and K.~Singh, ``Online control
  with adversarial disturbances,'' in \emph{36th International Conference on
  Machine Learning, ICML 2019}.\hskip 1em plus 0.5em minus 0.4em\relax
  International Machine Learning Society (IMLS), 2019, pp. 154--165.

\bibitem{agarwal2019logarithmic}
N.~Agarwal, E.~Hazan, and K.~Singh, ``Logarithmic regret for online control,''
  in \emph{Advances in Neural Information Processing Systems}, 2019, pp.
  10\,175--10\,184.

\bibitem{simchowitz2020improper}
M.~Simchowitz, K.~Singh, and E.~Hazan, ``Improper learning for non-stochastic
  control,'' \emph{arXiv preprint arXiv:2001.09254}, 2020.

\bibitem{hazan2020nonstochastic}
E.~Hazan, S.~Kakade, and K.~Singh, ``The nonstochastic control problem,'' in
  \emph{Algorithmic Learning Theory}, 2020, pp. 408--421.

\bibitem{lale2020logarithmic}
S.~Lale, K.~Azizzadenesheli, B.~Hassibi, and A.~Anandkumar, ``Logarithmic
  regret bound in partially observable linear dynamical systems,'' \emph{arXiv
  preprint arXiv:2003.11227}, 2020.

\bibitem{liu2008monte}
J.~S. Liu, \emph{Monte Carlo strategies in scientific computing}.\hskip 1em
  plus 0.5em minus 0.4em\relax Springer Science \& Business Media, 2008.

\bibitem{shahrampour2018distributed}
S.~Shahrampour and A.~Jadbabaie, ``Distributed online optimization in dynamic
  environments using mirror descent,'' \emph{IEEE Transactions on Automatic
  Control}, vol.~63, no.~3, pp. 714--725, 2018.

\end{thebibliography}


\end{document}